\ifx\documentstyle\undefined\else \fi 
\xdef\fmtversion{\fmtversion+CWEB3.1}

\let\:=\. 

\parskip 0pt 
\parindent 1em 

\let\mc=\ninerm 
\def\CEE/{{\mc C\spacefactor1000}}
\def\UNIX/{{\mc U\kern-.05emNIX\spacefactor1000}}
\def\TEX/{\TeX}
\def\CPLUSPLUS/{{\mc C\PP\spacefactor1000}}
\def\9#1{}
\font\eightrm=cmr8
\let\sc=\eightrm 
\let\mainfont=\tenrm
\let\cmntfont\tenrm
\font\tentex=cmtex10 
\fontdimen7\tentex=0pt 

\def\\#1{\leavevmode\hbox{\it#1\/\kern.05em}} 
\def\|#1{\leavevmode\hbox{$#1$}} 
\def\&#1{\leavevmode\hbox{\bf
  \def\_{\kern.04em\vbox{\hrule width.3em height .6pt}\kern.08em}%
  #1\/\kern.05em}} 
\def\.#1{\leavevmode\hbox{\tentex 
  \let\\=\BS 
  \let\{=\LB 
  \let\}=\RB 
  \let\~=\TL 
  \let\ =\SP 
  \let\_=\UL 
  \let\&=\AM 
  \let\^=\CF 
  #1\kern.05em}}
\def\){\discretionary{\hbox{\tentex\BS}}{}{}}
\def\ATL{\par\noindent\bgroup\catcode`\_=12 \postATL} 
\def\postATL#1 #2 {\bf letter \\{\uppercase{\char"#1}}
   tangles as \tentex "#2"\egroup\par}
\def\noATL#1 #2 {}
\def\noatl{\let\ATL=\noATL} 
\def\ATH{\X\kern-.5em:Preprocessor definitions\X}
\let\PB=\relax 

\chardef\AM=`\& 
\chardef\BS=`\\ 
\chardef\LB=`\{ 
\chardef\RB=`\} 
\def\SP{{\tt\char`\ }} 
\chardef\TL=`\~ 
\chardef\UL=`\_ 
\chardef\CF=`\^ 

\newbox\PPbox 
\setbox\PPbox=\hbox{\kern.5pt\raise1pt\hbox{\sevenrm+\kern-1pt+}\kern.5pt}
\def\PP{\copy\PPbox}
\newbox\MMbox \setbox\MMbox=\hbox{\kern.5pt\raise1pt\hbox{\sevensy\char0
 \kern-1pt\char0}\kern.5pt}
\def\MM{\copy\MMbox}
\newbox\MGbox 
\setbox\MGbox=\hbox{\kern-2pt\lower3pt\hbox{\teni\char'176}\kern1pt}
\def\MG{\copy\MGbox}
\def\MRL#1{\mathrel{\let\K==#1}}

\let\NULL=\Lambda
\mathchardef\AND="2026 
\newbox\MODbox \setbox\MODbox=\hbox{\eightrm\%}

\def\DC{\kern.1em{::}\kern.1em} 

\newbox\bak \setbox\bak=\hbox to -1em{} 
\newbox\bakk\setbox\bakk=\hbox to -2em{} 

\newcount\ind 
\def\1{\global\advance\ind by1\hangindent\ind em} 
\def\2{\global\advance\ind by-1} 
\def\3#1{\hfil\penalty#10\hfilneg} 
\def\4{\copy\bak} 
\def\5{\hfil\penalty-1\hfilneg\kern2.5em\copy\bakk\ignorespaces}
\def\6{\ifmmode\else\par 
  \hangindent\ind em\noindent\kern\ind em\copy\bakk\ignorespaces\fi}
\def\7{\Y\6} 
\def\8{\hskip-\ind em\hskip 2em} 

\newcount\gdepth 
\newcount\secpagedepth
\secpagedepth=3 
\newtoks\gtitle 
\newskip\intersecskip \intersecskip=12pt minus 3pt 
\let\yskip=\smallskip
\def\?{\mathrel?}
\def\note#1#2.{\Y\noindent{\hangindent2em\baselineskip10pt\eightrm#1~#2.\par}}
\def\lapstar{\rlap{*}}
\def\stsec{\rightskip=0pt 
  \sfcode`;=1500 \pretolerance 200 \hyphenpenalty 50 \exhyphenpenalty 50
  \noindent{\let\*=\lapstar\bf\secstar.\quad}}
\let\startsection=\stsec
\def\defin#1{\global\advance\ind by 2 \1\&{#1 } } 
\def\A{\note{See also section}} 
\def\B{\rightskip=0pt plus 100pt minus 10pt 
  \sfcode`;=3000
  \pretolerance 10000
  \hyphenpenalty 1000 
  \exhyphenpenalty 10000
  \global\ind=2 \1\ \unskip}
\def\C#1{\5\5\quad$/\ast\,${\cmntfont #1}$\,\ast/$}
\def\D{\defin{\#define}} 
\let\E=\equiv 
\let\G=\ge 
\let\I=\ne 
\let\K== 
\outer\def\M#1{\MN{#1}\ifon\vfil\penalty-100\vfilneg 
  \vskip\intersecskip\startsection\ignorespaces}
\outer\def\N#1#2#3.{\gdepth=#1\gtitle={#3}\MN{#2}
  \ifon\ifnum#1<\secpagedepth \vfil\eject 
    \else\vfil\penalty-100\vfilneg\vskip\intersecskip\fi\fi
  \message{*\secno} 
  \edef\next{\write\cont{\ZZ{#3}{#1}{\secno}
                   {\noexpand\the\pageno}}}\next 
  \ifon\startsection{\bf#3.\quad}\ignorespaces}
\def\MN#1{\par 
  {\xdef\secstar{#1}\let\*=\empty\xdef\secno{#1}}
  \ifx\secno\secstar \onmaybe \else\ontrue \fi
  \mark{{{\tensy x}\secno}{\the\gdepth}{\the\gtitle}}}
\def\Q{\note{This code is cited in section}} 
\def\Qs{\note{This code is cited in sections}} 
\let\R=\lnot 
\def\T#1{\leavevmode 
  \hbox{$\def\?{\kern.2em}%
    \def\$##1{\egroup_{\,\rm##1}\bgroup}
    \def\_{\cdot 10^{\aftergroup}}
    \let\~=\oct \let\^=\hex {#1}$}}
\def\U{\note{This code is used in section}} 
\def\Us{\note{This code is used in sections}} 
\let\W=\land 
\def\X#1:#2\X{\ifmmode\gdef\XX{\null$\null}\else\gdef\XX{}\fi 
  \XX$\langle\,${#2\eightrm\kern.5em#1}$\,\rangle$\XX}
\def\Y{\par\yskip}
\let\Z=\le
\let\ZZ=\let 
\let\*=*

\def\oct{\hbox{$^\circ$\kern-.1em\it\aftergroup\?\aftergroup}}
\def\hex{\hbox{$^{\scriptscriptstyle\#}$\tt\aftergroup}} 
\def\vb#1{\leavevmode\hbox{\kern2pt\vrule\vtop{\vbox{\hrule
        \hbox{\strut\kern2pt\.{#1}\kern2pt}}
      \hrule}\vrule\kern2pt}} 

\def\onmaybe{\let\ifon=\maybe} \let\maybe=\iftrue
\newif\ifon \newif\iftitle \newif\ifpagesaved

\def\lheader{\mainfont\the\pageno\eightrm\qquad\grouptitle\hfill\title\qquad
  \mainfont\topsecno} 
\def\rheader{\mainfont\topsecno\eightrm\qquad\title\hfill\grouptitle
  \qquad\mainfont\the\pageno} 
\def\grouptitle{\let\i=I\let\j=J\uppercase\expandafter{\expandafter
                        \takethree\topmark}}
\def\topsecno{\expandafter\takeone\topmark}
\def\takeone#1#2#3{#1}

\def\takethree#1#2#3{#3}
\def\nullsec{\eightrm\kern-2em} 

\let\page=\pagebody \raggedbottom
\def\normaloutput#1#2#3{\ifodd\pageno\hoffset=\pageshift\fi
 \shipout\vbox{
  \vbox to\fullpageheight{
  \iftitle\global\titlefalse
  \else\hbox to\pagewidth{\vbox to10pt{}\ifodd\pageno #3\else#2\fi}\fi
  \vfill#1}} 
  \global\advance\pageno by1}

\gtitle={\.{CWEB} output} 
\mark{\noexpand\nullsec0{\the\gtitle}}
\def\title{\expandafter\uppercase\expandafter{\jobname}}
\def\topofcontents{\centerline{\titlefont\title}\vskip.7in
  \vfill} 
\def\botofcontents{\vfill
  \centerline{\covernote}} 
\def\covernote{}
\def\contentspagenumber{0} 
\newdimen\pagewidth \pagewidth=6.5in 
\newdimen\pageheight \pageheight=8.7in 
\newdimen\fullpageheight \fullpageheight=9in 
\newdimen\pageshift \pageshift=0in 

\def\setpage{\hsize\pagewidth\vsize\pageheight} 
\def\contentsfile{\jobname.toc} 
\def\readcontents{\input \contentsfile}
\def\readindex{\input \jobname.idx}
\def\readsections{\input \jobname.scn}

\newwrite\cont
\output{\setbox0=\page 
  \openout\cont=\contentsfile
       \write\cont{\catcode `\noexpand\@=11\relax}   
  \global\output{\normaloutput\page\lheader\rheader}}
\setpage
\vbox to \vsize{} 

\def\ch{\note{The following sections were changed by the change file:}
  \let\*=\relax}
\newbox\sbox 
\newbox\lbox 
\def\inx{\par\vskip6pt plus 1fil 
  \def\page{\box255 } \normalbottom
  \write\cont{} 
       \write\cont{\catcode `\noexpand\@=12\relax}   
  \closeout\cont 
  \output{\ifpagesaved\normaloutput{\box\sbox}\lheader\rheader\fi
    \global\setbox\sbox=\page \global\pagesavedtrue}
  \pagesavedfalse \eject 
  \setbox\sbox\vbox{\unvbox\sbox} 
  \vsize=\pageheight \advance\vsize by -\ht\sbox 
  \hsize=.5\pagewidth \advance\hsize by -10pt
  \parfillskip 0pt plus .6\hsize 
  \def\lr{L} 
  \output{\if L\lr\global\setbox\lbox=\page \gdef\lr{R}
    \else\normaloutput{\vbox to\pageheight{\box\sbox\vss
        \hbox to\pagewidth{\box\lbox\hfil\page}}}\lheader\rheader
    \global\vsize\pageheight\gdef\lr{L}\global\pagesavedfalse\fi}
  \message{Index:}
  \parskip 0pt plus .5pt
  \outer\def\I##1, {\par\hangindent2em\noindent##1:\kern1em} 
  \def\[##1]{$\underline{##1}$} 
  \rm \rightskip0pt plus 2.5em \tolerance 10000 \let\*=\lapstar
  \hyphenpenalty 10000 \parindent0pt
  \readindex}
\def\fin{\par\vfill\eject 
  \ifpagesaved\null\vfill\eject\fi 
  \if L\lr\else\null\vfill\eject\fi 
  \parfillskip 0pt plus 1fil
  \def\grouptitle{NAMES OF THE SECTIONS}
  \let\topsecno=\nullsec
  \message{Section names:}
  \output={\normaloutput\page\lheader\rheader}
  \setpage
  \def\note##1##2.{\quad{\eightrm##1~##2.}}
  \def\Q{\note{Cited in section}} 
  \def\Qs{\note{Cited in sections}} 
  \def\U{\note{Used in section}} 
  \def\Us{\note{Used in sections}} 
  \def\I{\par\hangindent 2em}\let\*=*
  \readsections}
\def\con{\par\vfill\eject 
  \rightskip 0pt \hyphenpenalty 50 \tolerance 200
  \setpage \output={\normaloutput\page\lheader\rheader}
  \titletrue 
  \pageno=\contentspagenumber
  \def\grouptitle{TABLE OF CONTENTS}
  \message{Table of contents:}
  \topofcontents
  \line{\hfil Section\hbox to3em{\hss Page}}
  \let\ZZ=\contentsline
  \readcontents\relax 
  \botofcontents \end} 
\def\contentsline#1#2#3#4{\ifnum#2=0 \smallbreak\fi
    \line{\consetup{#2}#1
      \rm\leaders\hbox to .5em{.\hfil}\hfil\ #3\hbox to3em{\hss#4}}}
\def\consetup#1{\ifcase#1 \bf 
  \or 
  \or \hskip2em 
  \or \hskip4em 
  \or \hskip6em 
  \or \hskip8em 
  \or \hskip10em 
  \else \hskip12em \fi} 
\def\noinx{\let\inx=\end} 
\def\nosecs{\let\FIN=\fin \def\fin{\let\parfillskip=\end \FIN}}
\def\nocon{\let\con=\end} 
\def\today{\ifcase\month\or
  January\or February\or March\or April\or May\or June\or
  July\or August\or September\or October\or November\or December\fi
  \space\number\day, \number\year}
\newcount\twodigits
\def\hours{\twodigits=\time \divide\twodigits by 60 \printtwodigits
  \multiply\twodigits by-60 \advance\twodigits by\time :\printtwodigits}
\def\gobbleone1{}
\def\printtwodigits{\advance\twodigits100
  \expandafter\gobbleone\number\twodigits
  \advance\twodigits-100 }
\def\TeX{{\ifmmode\it\fi
   \leavevmode\hbox{T\kern-.1667em\lower.424ex\hbox{E}\hskip-.125em X}}}
\def\,{\relax\ifmmode\mskip\thinmuskip\else\thinspace\fi}
\def\datethis{\def\startsection{\leftline{\sc\today\ at \hours}\bigskip
  \let\startsection=\stsec\stsec}}


\def\topofcontents{
 \null\vskip-35pt
 \line{\eightrm Journal of Experimental Algorithmics\hfill}
 \vfill
 \vfill
 \centerline{\bf Irredundant Intervals}
 \bigskip
 \centerline{Donald E. Knuth}
 \smallskip
 \centerline{Computer Science Department, Stanford University}
 \vskip 4\bigskipamount
 \centerline{\vtop{\hsize=.8\hsize
    \noindent{\bf Abstract.}\quad
    This expository note presents simplifications of a theorem
    due to Gy\H{o}ri and an algorithm due to Franzblau and Kleitman:
    Given a family $F$ of $m$ intervals on a linearly ordered set
    of $n$ elements, we can construct in $O(m+n)^2$ steps an
    irredundant subfamily having maximum cardinality, as well as
    a generating family having minimum cardinality. The algorithm
    is of special interest because it solves a problem analogous
    to finding a maximum independent set, but on a class of objects
    that is more general than a matroid.
    This note is also a complete, runnable computer program,
    which can be used for experiments in conjunction with the
    public-domain software of {\sl The Stanford GraphBase}.}}
 \vskip.7in
 \vfill} 

\let\oldcon=\con
\let\con=\bye

\newread\tocfile \openin\tocfile=\contentsfile\relax
\ifeof\tocfile
 \message{Please run this file again thru TeX to get the table of contents!}
\else \pageno=0 \titletrue {\let\end=\relax \oldcon} \fi

\def\lheader{\mainfont\the\pageno\eightrm\qquad\grouptitle\hfill
  Journal of Experimental Algorithmics\hfill
  DONALD E. KNUTH\qquad
  \mainfont\topsecno} 
\def\rheader{\mainfont\topsecno\eightrm\qquad IRREDUNDANT INTERVALS\hfill
  Journal of Experimental Algorithmics\hfill\grouptitle
  \qquad\mainfont\the\pageno} 

\font\slttten=cmsltt10


\def\CR{1} 
\def\FJ{2} 
\def\FK{3} 
\def\Gy{4} 
\def\KB{5} 
\def\vone{6} 
\def\vthree{7} 
\def\sgb{8} 
\def\cweb{9} 
\def\Lu{10} 


\def\[#1,#2){[#1\,\mathpunct{\ldotp\ldotp}#2)} 
\def\:{{\vert\mskip1mu}} 
\def\]{{\downarrow}\mskip1mu} 
\let\seq=\subseteq 
\let\union=\bigcup
\def\(#1){^{\mskip1mu(#1)}}

\let\int=\[ 
\def\proclaim #1. #2\par{{\bf #1.}\enspace{\sl#2}\par}
\def\Proof{\smallskip\noindent Proof:\enspace}
\def\callthis#1{\edef#1{\number\secno}}


\newdimen\unit \unit=1.5em
\def\intpic #1[#2..#3)#4\\{\hbox{\kern#2\unit\llap{$#1$\ \ }%
 \hbox to0pt{\hss$\bullet$\hss}%
 \dimen0=#3\unit \advance\dimen0-#2\unit
 \hbox to\dimen0{\cleaders\hbox{$\mkern-2mu\smash-\mkern-2mu$}\hfil
   $\mkern-7mu \smash- \mkern1.5mu$}%
 \hbox to0pt{\hss$\circ$\hss} \ $#4$}}

\N{1}{1}Introduction. Let's say that a family of sets is {\it irredundant\/} if
its members can be arranged in a sequence with the following property:
Each set contains a point that isn't in any of the preceding sets.

If $F$ is a family of sets, we write $F^\cup$ for the family of all nonempty
unions of elements of $F$. When $F$ and $G$ are families with
$F\subseteq G^\cup$, we say that $G$ {\it generates\/} $F$. If $F$ is
irredundant and $G$ generates $F$, we obviously have $\vert F\vert\le
\vert G\vert$, because each set in the sequence requires a new generator.

In the special case that the members of $F$ are intervals of the real line,
Andr\'as Frank
conjectured that the largest irredundant subfamily of $F$ has the same
cardinality as $F$'s smallest generating family. This conjecture was
proved by Ervin Gy\H{o}ri~[\Gy], who noted that such a result was a minimax
theorem of a new type, apparently unrelated to any of the other well-known
minimax theorems of graph theory and combinatorics. A constructive proof
was found shortly afterwards by Franzblau and Kleitman~[\FK], who
sketched an algorithm to find a generating family and irredundant subfamily of
equal cardinality.
(Gy\H{o}ri, Franzblau, and Kleitman were led to these results while studying
the more general problem of finding a minimum number of subrectangles that
cover a given polygon. Further information about polygon covers appears
in [\FK] and~[\CR].)

The purpose of this note is to describe the beautiful algorithm of
Franzblau and Kleitman in full detail.
Indeed, the \.{CWEB} source file that generated this document is a computer
program that can be used in connection with the Stanford GraphBase~[\sgb] to
find maximum irredundant subfamilies and minimum generating families of any
given collection of intervals. Perhaps this new exposition will
shed new light on the class of optimization problems
for which an efficient algorithm exists.

According to the conventions of \.{CWEB}~[\cweb], the sections of this
document are
sequentially numbered 1, 2, 3, etc. In this respect we are returning to a
style of exposition used by Euler and Gauss and their contemporaries.
A~\.{CWEB} program is also essentially a hypertext; therefore this document
may also be regarded as experimental in another sense, as an attempt
to find new forms of exposition appropriate to modern technology.

Note: Gy\H{o}ri used the term ``U-increasing'' for an irredundant
family; Franzblau and Kleitman called such intervals ``independent.''
Since a family of sets is a hypergraph, it seems unwise to deviate
from the standard meaning of independent edges, yet ``U-increasing'' is
not an especially appealing alternative. We will see momentarily that
the term ``irredundant'' is quite natural in theory and practice.

\fi

\M{2}A far-reaching generalization of Gy\H{o}ri's theorem was proved recently
by
Frank and Jord\'an~[\FJ], who introduced a large new family of minimax
theorems related to linking systems. In particular, Frank and Jord\'an
extended Gy\H{o}ri's results to intervals on a circle instead of a line.
But no combinatorial algorithm is known as yet for the circular case.

\callthis\intro
Can any or all of the Franzblau/Kleitman methods be ``lifted'' to such more
general problems? We will return to this tantalizing question after becoming
familiar with Franzblau and Kleitman's remarkable algorithm.

\fi

\N{1}{3}Theory. It is wise to study the theory underlying the
Franzblau/Kleitman algorithm before getting into the program itself.

\fi

\M{4}A family of sets is called {\it redundant\/} if it is not irredundant.
Any family that contains a redundant subfamily is redundant, since any
family contained in an irredundant family is irredundant.

\fi

\M{5}If $F$ is a family of sets and $s$ is an arbitrary set, let $F\:s$ denote
the
sets of $F$ that are contained in~$s$. This operation is left-associative
by convention: $F\:s\:t=(F\:s)\:t=F\:(s\cap t)$.

We also write $\union F$ for $\union\{f\mid f\in F\}$; thus
$F\:\,{\union F}=F$.

(An index to all the main notations and definitions that we will use
appears at the end of this note.)

\fi

\M{6}\proclaim Lemma. A finite family $F$ is redundant if and only if there is
a nonempty set $s$ such that every point of~$s$ belongs to at least two members
of $F\:s$. {\rm(The set $s$ need not belong to~$F$.)}

\Proof If $F$ is irredundant there is no such $s$, because $F\:s$ is
irredundant and its last set in the assumed sequence contains a
point that isn't in any of the others. But if $F$ is redundant, it contains
a minimal redundant subfamily $F_0$; then we have
$$\textstyle f\seq \union\,(F_0\setminus\{f\})
\qquad\hbox{for all $f\in F_0$}\,,$$
since $F_0\setminus\{f\}$ is irredundant. It follows that every point of
$s=\union F_0$ is contained in at least two members of $F_0$, hence in
at least two members of $F\:s$ (since $F_0\seq F\:s$).

\fi

\M{7}\proclaim Corollary. A finite family $F$ of intervals on a line
is redundant if and
only if there is an interval $s$ such that every point of $s$ belongs to
at least two intervals of $F\:s$. {\rm(The set $s$ need not belong to~$F$.)}

\Proof Intervals are nonempty.
By the proof of the preceding lemma, it suffices to consider sets~$s$
that can be written $\union F_0$ for some minimal redundant subfamily~$F_0$.
In the special case of intervals, $\union F_0$ must be a single interval;
otherwise $F_0$ would not be minimal.

\fi

\M{8}Henceforth we will restrict consideration to finite families $F$ of
intervals on a linearly ordered set.
It suffices, in fact, to deal with integer elements; we
will consider subintervals of the $n$-element set $\[0,n)$. (The
notation $\[a,b)$ stands here for the set of all integers $x$ such that
$a\le x<b$.)

If $F$ is a family of sets and $x$ is a point, we will write $N_x\,F$ for the
number of sets that contain $x$. The corollary just proved can therefore be
stated as follows: ``$F$ is irredundant if and only if
every interval $s\seq\union F$ contains a point $x$ with $N_x\,F\:s\le1$.''
This characterization provides a polynomial-time test for irredundancy.

\fi

\M{9}\callthis\bst
Irredundant intervals have an interesting connection to the familiar
computer-science concept of {\it binary search trees\/} (see, for example,
[\vthree, \S6.2.2]): A family of intervals is irredundant if and only if
we can associate its intervals with a binary tree whose
nodes are each labeled with an integer~$x$ and an interval containing~$x$.
All nodes in the left subtree of such a node correspond to intervals
that are strictly less than~$x$, in the sense that all elements of those
intervals are $<x$; all nodes in the right subtree correspond
to intervals that are strictly greater than~$x$. The root of the binary tree
corresponds to the interval that is last in the assumed irredundant
ordering. Its distinguished integer~$x$ is an element that appears in
no other interval.

Given such a tree, we obtain a suitable irredundant ordering by traversing it
recursively from the leaves to the root, in postorder [\vone, \S2.3.1].
Conversely, given an irredundant ordering, we can construct a binary tree
recursively, proceeding from the root to the leaves.

\fi

\M{10}An example might be helpful at this point. Suppose $n=9$ and
$$f_1=\[0,8)\,,\quad f_2=\[0,7)\,,\quad f_3=\[1,6)\,,\quad
f_4=\[1,5)\,,\quad f_5=\[3,9)\,,\quad f_6=\[2,9)\,.$$
Then $\{f_1,f_2,f_3,f_4,f_5\}$ and $\{f_1,f_3,f_5,f_6\}$ are irredundant.
(Indeed, a family of intervals is irredundant whenever its members have no
repeated left endpoints or no repeated right endpoints.) These subfamilies are
in fact {\it maximally\/} irredundant---they become redundant when any
other interval of the family is added. Therefore maximal irredundant
subfamilies need not have the same cardinality; irredundant subfamilies do not
form the independent sets of a matroid.

On the other hand, irredundant sets of intervals do have matroid-like
properties. For example, if $F$ is irredundant and $F\cup\{g\}$ is redundant,
there is an $f\ne g$ such that $F\cup \{g\}\setminus\{f\}$ is irredundant.
(The proof is by induction on $\vert F\vert$: There is
$x\in f\in F$ such that $F=F_l\cup\{f\}\cup F_r$, where $F_l$ and $F_r$
correspond to the left and right subtrees of the root in the binary tree
representation. If $x\in g$, the family $F_l\cup\{g\}\cup F_r$ is irredundant;
if $x<g$, there is $f'\in F_r$ with
$F_l\cup\{f\}\cup (F_r\cup\{g\}\setminus\{f'\})$
irredundant, by induction; if $x>g$ there is $f'\in F_l$ with
$(F_l\cup\{g\}\setminus\{f'\})\cup\{f\}\cup F_r$ irredundant.)
Such near-matroid behavior makes families of intervals especially instructive.

\fi

\M{11}Let's say that an interval $s$ is {\it good\/} for $F$ if $N_x\,F\:s\le1$
for some $x\in s$; otherwise $s$ is {\it bad}. Franzblau and Kleitman
introduced a basic reduction procedure
for any family $F$ of intervals that possesses a bad interval~$s$. Their
procedure is analogous to modification along an augmenting path in other
combinatorial algorithms.

Let $\[a_1,b_1)$, \dots, $\[a_k,b_k)$ be the maximal intervals in~$F\:s$,
ordered so that $a_1<\cdots<a_k$ and $b_1<\cdots<b_k$. (Notice that
$s=\[a_1,b_k)$.) For example, we might have the following picture:
$$\vbox{
\halign{\kern5\unit\hfil#\hfil\cr
$s$\cr
\noalign{\vskip-.4\baselineskip}
\hbox to 20\unit{\downbracefill}\cr}
\medskip
\intpic a_1[5..9)b_1\\
\intpic a_2[7..15)b_2\\
\intpic a_3[9..19)b_3\\
\intpic a_4[14..21)b_4\\
\intpic a_k[17..25)b_k\\
\intpic [5..8)\hbox to 12\unit{\hss
(non-maximal intervals make $s$ bad)\ \ \ \ \hss}\\
\vskip-\baselineskip \intpic [20..25)\\
\hbox to 30\unit{\hfil}\vskip-\baselineskip 
}$$
If $a_{j+1}<b_j$ for $1\le j<k$, the family $F$ {\it reduced in\/} $s$ is
defined to be
$$F\]s \;=\; F \,\setminus\, \{\[a_1,b_1),\ldots,\[a_k,b_k)\}
\,\cup\,      \{\[a_2,b_1),\ldots,\[a_k,b_{k-1})\}\,.$$
In the simplest case we have $k=1$ and the
reduced family is simply $F\setminus\{\[a_1,b_1)\}$.

If $s$ is a {\it minimal\/} bad interval for $F$, we can prove that
$a_{j+1}<b_j$ for $1\le j<k$, hence $F\]s$ is well-defined. Indeed,
point~$a_{j+1}$ must be in some interval $\[c,d)$ other than
$\[a_{j+1},b_{j+1})$, since $s$ is bad. We can assume that $c<a_{j+1}$;
otherwise all intervals of $F\:s$ would be contained in $\[a_1,a_{j+1})$ or
$\[a_{j+1},a_k)$, and both of these subintervals would be bad, contradicting
the minimality of~$s$. But $c<a_{j+1}$ implies that $\[c,d)$ is contained in
some maximal $\[a_i,b_i)$ with $i\le j$. Hence $a_{j+1}<b_i\le b_j$.

The notation $F\]s$ is defined to be left-associative, like $F\:s$; that is,
$F\]s\]t=(F\]s)\]t$ and $F\]s\,\:t=(F\]s)\:t$.

\fi

\M{12}\callthis\lemzero
\proclaim Lemma. If $s$ is a minimal bad interval for $F$, we have
$F\seq(F\]s)^\cup$.

\Proof We must show that $\[a_j,b_j)\in(F\]s)^\cup$ for $1\le j\le k$. Let
$G=F\]s\,\:\[a_j,b_j)$; we will prove that $\[a_j,b_j)=\union G$. If
$x\in\[a_j,b_j)$, the badness of~$s$ implies that $x\in t$ for some $t\in F\:s$
with $t\ne\[a_j,b_j)$; let $t$ be contained in the maximal interval
$\[a_i,b_i)$.
If $i=j$, we have $t\in G$;
if $i<j$, we have $x\in\[a_j,b_i)\seq\[a_j,b_{j-1})\in G$; and
if $i>j$, we have $x\in\[a_i,b_j)\seq\[a_{j+1},b_j)\in G$.

\fi

\M{13}\callthis\lemone \proclaim Lemma.
Suppose $s$ is a minimal bad interval for $F$, while
$t$ is a good interval. Then $t$ is good also for $F\]s$.

\Proof Let $s=\[a,b)$, $t=\[c,d)$, $l=\min(a,c)$, $r=\max(b,d)$.
Suppose $N_x\,F\:\[l,d)\ge2$ and $N_x\,F\:\[c,r)\ge2$ for all
$x\in t$.  Then $l=a<c<d <b=r$, because $t$ is good for $F$. By the
minimality of $s$, there is some $x\in\[a,d)$ with $N_x\,F\:\[a,d)\cdot
x\le1$. Since $N_x\,F\:\[a,d)\ge2$ we have $x<c$. Furthermore,
$x$ is in some interval of $(F\:s)\setminus F\:\[a,d)$,
because $s$ is bad; so $x$ is in some maximal $\[a_j,b_j)$ with
$a\le a_j\le x<c<d<b_j\le b$. It follows that none of the intervals
$\[a_1,b_1)$, \dots, $\[a_k,b_k)$, $\[a_2,b_1)$, \dots, $\[a_k,b_{k-1})$
are contained in~$t$, hence $F\]s\,\:t=F\:t$.

On the other hand, suppose $N_x\,F\:\[l,d)\le1$ for some $x\in t$. We will show
that $N_x\,F\]s\,\:t\le1$. This assertion can fail only if $x$ lies in some
interval $\[a_{j+1},b_j)\seq t$ newly added to $F\]s$.
Then $x\in\[a_j,b_j)\seq\[l,d)$, and $x\in\[a_{j+1},b_{j+1})\notin\[l,d)$,
hence $b_j\le d<b_{j+1}$; it follows that $j$ is uniquely determined, and
the only interval containing $x$ in $F\]s\,\:t$ is $\[a_{j+1},b_j)$.
A similar argument applies if $N_x\,F\:\[c,r)\le1$.

\fi

\M{14}\callthis\pretrick \proclaim Corollary.
If $s$ is a minimal bad interval for $F$, we have $N_x\,F\]s=N_x\,F-N_x\,s$.

\Proof The proof of the preceding lemma shows in particular that none of the
intervals $\[a_{j+1},b_j)$ are already present in~$F$ before the reduction.
And if $x\in s$, suppose $x$~lies in $\[a_i,b_i)$, $\[a_{i+1},b_{i+1})$, \dots,
$\[a_j,b_j)$; then it lies in $\[a_{i+1},b_i)$, \dots, $\[a_j,b_{j-1})$ after
reduction, a net change of~$-1$.

\goodbreak

\fi

\M{15}The Franzblau/Kleitman algorithm has a very simple outline: We let
$G_0=F$ and repeatedly set $G_{k+1}=G_k\]s_k$, where $s_k$ is the leftmost
minimal bad interval for $G_k$, until we finally reach a family $G_r$ in which
no bad intervals remain. This must happen sooner or later, because $\vert
G_k\vert =\vert F\vert-k$. The final irredundant family $G=G_r$ generates~$F$,
because $F\seq G_k^\cup$ for all~$k$ by the lemma of \S\lemzero.
Franzblau and Kleitman proved the nontrivial
fact that $\vert G\vert$ is the size of the maximum irredundant subfamily
of~$F$; hence $G$ is a minimum generating family. \callthis\outline

\fi

\M{16}It is tempting to try to prove the optimality of $G$ by a simpler,
inductive approach in which we ``grow'' $F$ one interval at a time,
updating its maximum irredundant set and minimum generating set appropriately.
But experiments show that the maximum irredundant set can change drastically
when $F$ receives a single new interval, so this direct primal-dual approach
seems doomed to failure. The indirect approach is more difficult to prove,
but no more difficult to program. So we will proceed to develop further
properties of Franzblau and Kleitman's reduction procedure~[\FK]. The key fact
is a remarkable theorem that we turn to~next.

\fi

\M{17}\proclaim Theorem. The same final family $G=G_r$ is obtained when $s_k$
is
chosen to be an arbitrary (not necessarily leftmost) minimal bad interval
of~$G_k$ in the reduction algorithm. Moreover, the same multiset
$\{s_0,\ldots,s_{r-1}\}$ of minimal bad intervals arises, in some order,
regardless of the choices made at each step.

\Proof We use induction on $r$, the maximum number of steps to convergence
among all reduction procedures that begin with a family $F$.
If $r=0$, the result is trivial,
and if $F$ has only one minimal bad interval the result is immediate
by induction. Suppose therefore that $s$ and $t$ are distinct minimal
bad intervals of $F$. We will prove later that $t$ is a minimal bad interval
for $F\]s$, and that $F\]s\]t=F\]t\]s$. Let $r'$ be the maximum distance
to convergence from $F\]s$, and $r''$ the maximum from $F\]t$; then $r'$
and $r''$ are less than $r$, and induction proves that the final result
from $F\]s$ is the final result from $F\]s\]t=F\]t\]s$, which is the final
result from $F\]t$. (Readers familiar with other reduction algorithms, like
that of [\KB], will recognize this as a familiar ``diamond lemma'' argument.
We construct a diamond-shaped diagram with four vertices: $F$, $F\]s$, $F\]t$,
and a common outcome of $F\]s$ and $F\]t$.)
This completes the proof, except for two lemmas that will be demonstrated
below; their proofs have been deferred so that we could motivate them first.

\fi

\M{18}This theorem and the lemma of \S\lemone\ have an important
corollary: {\sl Let $S=\{s_0,\ldots,s_{r-1}\}$ be the multiset of minimal
bad intervals determined by the algorithm from~$F$, and let $t$ be any
interval. Then $S\:t$ is the multiset of minimal bad intervals determined by
the algorithm from~$F\:t$.} This holds
because an interval $s\seq t$ is bad for $F$ if and only if it is bad for
$F\:t$. Minimal bad intervals within $t$ never appear again once they are
removed, and we can remove them first.

Reducing a minimal bad interval $s$ when $s$ is contained in a bad
interval~$t$ may make $t$ good, or leave it bad, or make it minimally bad.
If $s$ is minimally bad for~$F$, it might also be minimally bad for $F\]s$.

\fi

\M{19}Now we are ready for the {\it coup de gr\^ace\/} and the {\it pi\`ece
de r\'esistance\/}:
After the reduction algorithm has computed the irredundant generating
family~$G=G_r$ and the multiset $S$ of minimal bad intervals, we
can construct an irredundant subfamily $F'$ of $F$ with $\vert F'\vert=
\vert G\vert$ by constructing a binary search tree as described in \S\bst.
The procedure is recursive, starting with an initial interval $t=\[0,n)$
that contains~$F$: The tree defined for $F\:t$ is empty if
$F\:t$ is empty. Otherwise it has a root node labeled with $x$
and with any interval of $F\:t$ containing~$x$, where $x$ is an integer such
that $N_x\,G^{(t)}=1$;
here $G^{(t)}$ is the final generating set that is obtained
when the reduction procedure is applied to $F\:t$. A suitable interval
containing~$x$ exists, because every element of $G^{(t)}$ is an intersection
of intervals in $F\:t$. The left subtree of the root node is
the binary search tree for $F\:\bigl(t\cap\[0,x)\bigr)$; the
right subtree is the binary search tree for $F\:\bigl(t\cap\[x+1,n)\bigr)$.

The number of nodes in this tree is $\vert G\vert$. For if $x$ is the integer
in the label of the root, $G$ has one interval containing $x$, and its
other intervals are $G\:\[0,x)$ and $G\:\[x+1,n)$. The family $G^{(t)}$ is not
the same as $G\:t$; but we do have $\vert G^{(t)}\vert=
\bigl\vert G\:t\bigr\vert$ when $t$ has the special form $\[0,x)$ or
$\[x+1,n)$, because in such cases
$F\]s\,\:t$ has the same cardinality as $F\:t$ when $s$ is a minimal bad
interval and $s\not\seq t$.
For example, if $t=\[0,x)$ and $b_j\le x<b_{j+1}$, we have $F\]s\,\:t=
F\:t\setminus\{\[a_1,b_1),\ldots,\[a_j,b_j)\}\cup\{\[a_2,b_1),\ldots,
\[a_{j+1},b_j)\}$.

\fi

\M{20}It is not necessary to compute each $G^{(t)}$ from scratch by starting
with
$F\:t$ and applying the reduction algorithm until it converges, because
the binary tree construction algorithm requires only a knowledge
of the incidence function $N_x\,G^{(t)}$. This function is easy to compute,
because $N_x\,F\]s=N_x\,F-N_x\,s$ by \S\pretrick; therefore
$$N_x\,G^{(t)}=N_x\,F\:t-N_x\,S\:t\,.$$ \callthis\trick

\fi

\M{21}All the basic ideas of Franzblau and Kleitman's algorithm have now
been explained. But we must still carry out a careful analysis of some
fine points of reduction that were claimed in the proof of the main theorem.
If $s$ and $t$ are distinct minimal bad intervals, the lemma of \S\lemone\
implies that no bad subintervals of~$t$ appear in $F\]s$; we also need
to verify that $t$ itself remains bad.

\smallskip
\noindent\proclaim Lemma. If $s$ is a minimal bad interval for $F$ and
$t$ is a bad interval such that $s\not\seq t$, then $t$ is bad for $F\]s$.

\Proof Let $s=\[a,b)$ and $t=\[c,d)$. We can assume by left-right symmetry
that $a<c$. Then $b<d$, by minimality of~$s$.
Assume that $t$ isn't bad for $F\]s$. The subfamily $F\:t$ must contain
at least one of the maximal intervals $\[a_j,b_j)$ of $F\:s$ that are
deleted during the reduction; hence $c\le a_j<b_{j-1}\le b$.

Let $j$ be minimal with $a_j\ge c$. Then
$$F\]s\,\:t=F\:t\setminus \{\[a_j,b_j),\ldots,\[a_k,b_k)\}
\cup \{\[a_j,b_{j-1}),\ldots,\[a_k,b_{k-1})\}\,;$$
so the elements of $t$ that are covered once less often are the elements of
$\[b_{j-1},b_k)$. Suppose $y\in\[b_{i-1},b_i)$ for some $i\ge j$. Then
$y\in\[a_i,b_i)\seq s\cap t$, and $y$ is in some other interval $f\seq s$.
The maximal interval containing $f$ must be $\[a_l,b_l)$ for some $l\ge i$,
hence $f\seq s\cap t$. Thus $N_y\,F\:(s\cap t)\ge 2$ for all $y\in
\[b_{j-1},b_k)$. But $s\cap t$ is good, so there must be a point $x\in
\[c,b_{j-1})$ with $N_x\,F\:(s\cap t)\le1$. We also have $N_x\,F\:t\ge 2$,
since $t$ is bad, so there's an interval in $F\:t\setminus F\:(s\cap t)$
that contains~$x$. This interval contains
$\[x,b_k)$. Consequently $N_y\,F\:t\ge3$ for all $y\in\[b_{j-1},b_k)$.

\fi

\M{22}\proclaim Lemma. If $s$ and $t$ are minimal bad intervals for $F$ and
$s\ne t$, we have $F\]s\]t=F\]t\]s$.

\Proof Let the maximal intervals in $F\:s$ and $F\:t$ be $\[a_1,b_1)$,
\dots, $\[a_k,b_k)$ and $\[c_1,d_1)$, \dots, $\[c_l,d_l)$, respectively,
where $a_1<c_1$. The lemma is obvious unless $F\:(s\cap t)$ is nonempty, so we
assume that $c_1<b_k<d_l$. Let $x\in s\cap t$ have $N_x\,F\:(s\cap t)=1$,
and let $f$ be the interval of $F\:(s\cap t)$ that contains~$x$. Let $p$
be maximal with $a_p<c_1$, and let $q$ be minimal with $d_q>b_k$. Since
$N_x\,F\:s>1$, there is an interval $\[a_j,b_j)$ containing~$x$
with $j\le p$; thus
$x\in\[a_p,b_p)$. Similarly $x\in\[c_q,d_q)$. Furthermore, if $p<k$ we
have $\[a_{p+1},b_{p+1})\seq s\cap t$; hence
either $\[a_{p+1},b_{p+1})=f$ or $a_{p+1}>x$. If $q>1$ we have either
$\[c_{q-1},d_{q-1})=f$ or $d_{q-1}\le x$.

If $p=k$ or
$f\ne\[a_{p+1},b_{p+1})$, any newly added intervals $\[a_{j+1},b_j)$ for
$p\le j<k$ in $F\]s$ are properly contained in $\[c_q,d_q)$, so they remain
in $F\]s\]t$. Thus we can easily describe the compound operation $F\]s\]t$ in
detail:
$$\displaylines{\hskip6em\hbox{Delete $\[a_1,b_1)$, \dots, $\[a_k,b_k)$;\quad
insert $\[a_2,b_1)$, \dots, $\[a_k,b_{k-1})$;}\hfill\cr
\hfill\hbox{delete $\[c_1,d_1)$, \dots, $\[c_l,d_l)$;\quad
insert $\[c_2,d_1)$, \dots, $\[c_l,d_{l-1})$.}\hskip6em\cr}$$
No two of these intervals are identical, so $F\]t\]s$ gives the same result.
(If $f=\[c_{q-1},d_{q-1})$, the family $F\]t$ has $f$ replaced by
$\[c_q,d_{q-1})\seq f\seq\[a_i,b_i)$ for some $i\le p$,
so $\[c_q,d_{q-1})$ is not maximal in $F\]t\,\:s$.)

The remaining case $f=\[a_{p+1},b_{p+1})=\[c_{q-1},d_{q-1})$ needs to be
considered specially, since we can't delete this interval twice.
The following picture might help clarify the situation:
$$\vbox{
\halign{\hfil#\hfil\cr
$s$\cr
\noalign{\vskip-.4\baselineskip}
\hbox to 19\unit{\downbracefill}\cr}
\medskip
\intpic a_1[0..5)b_1\\
\intpic a_2[2..12)b_2\\
\intpic a_p[3..15)b_p\\
\intpic (\;a_{p+1}[6..16)b_{p+1}\quad
\hbox{same as $\[c_{q-1},d_{q-1})$ below$\;)$}\\
\intpic a_k[11..19)b_k\\
\smallskip
\intpic c_1[4..7)d_1\\
\intpic c_{q-1}[6..16)d_{q-1}\\
\intpic c_q[8..23)d_q\\
\intpic c_l[14..26)d_l\\
\smallskip
\halign{\kern4\unit\hfil#\hfil\cr
\hbox to 22\unit{\upbracefill}\cr
$t$\cr}
\moveright 9\unit\hbox{\smash{\vbox to 12.5\baselineskip{
\cleaders\hbox to0pt{\hss\vrule height.5ex
\vrule height 1ex depth .5ex width0pt\hss}\vfil
\medskip\hbox to0pt{\hss$x$\hss}}}}
}$$

Suppose $g$ is an interval of $F\:t$ that is contained in $\[c_j,d_j)$ if and
only if $j=q-1$. Then $g$ contains a point $<c_q$. If $x\in g$ we have $g=f$,
since $g\seq s\cap t$. Otherwise $g\seq\[c_{q-1},x)\seq\[c_{q-1},b_p)=
\[a_{p+1},b_p)$. It follows that the maximal intervals of $F\]s\,\:t$ are
$$\[c_1,d_1),\;\ldots,\;\[c_{q-2},d_{q-2}),\;\[c_{q-1},b_p),\;
\[c_q,d_q),\;\ldots,\;\[c_l,d_l).$$
These intervals are replaced in $F\]s\]t$ by
$$\[c_2,d_1),\;\ldots,\;\[c_{q-1},d_{q-2}),\;\[c_q,b_p),\;
\[c_{q+1},d_q),\;\ldots,\;\[c_l,d_{l-1}).$$
Thus $F\]s\]t$ is formed almost as in the previous case, but with $\[a_{p+1},
b_p)$ and $\[c_q,d_{q-1})$ replaced by $\[c_q,b_p)$. And we get precisely
the same intervals in $F\]t\]s$.

(Is there a simpler proof?)

\fi

\N{1}{23}Practice. The computer program in the remainder of this note
operates on a family of intervals defined by a graph on $n+1$
vertices $\{0,1,\ldots,n\}$. We regard an edge between $u$ and $v$
as the half-open interval $\[u,v)$, when $u<v$.

Graphs are represented as in the algorithms of the Stanford GraphBase~[\sgb],
and the reader of this program is supposed to be familiar with the elementary
conventions of that system.

The program reads two command-line parameters,
$m$ and $n$, and an optional third parameter representing a random-number
seed. (The seed value is zero by default.) The Franzblau/Kleitman algorithm
is then applied to the graph \PB{$\\{random\_graph}(\|n+\T{1},\|m,\T{0},\T{0},%
\T{0},\T{0},\T{0},\T{0},\T{0},\\{seed})$},
a~random graph with vertices $\{0,1,\ldots,n\}$ and $m$ edges.
Alternatively, the user can specify an arbitrary graph as input by typing
the single command-line parameter \.{-g}$\langle\,$filename$\,\rangle$;
in this case the named file should describe the graph in \PB{\\{save\_graph}}
format
(as in the {\sc MILES\_SPAN} program of~[\sgb]).

When the computation is finished, a minimal generating family and a
maximal irredundant subfamily will be printed on the standard output file.

If a negative value is given for \PB{\|n}, the random graph is reversed from
left to right; each interval $\[a,b)$ is essentially replaced by $\[-b,-a)$
(but minus signs are suppressed in the output).
This feature lends credibility to the correctness of our
highly asymmetric algorithm and program, because we can verify the fact that
the minimum generating family of the mirror image of~\PB{\|F} is indeed
the mirror image of \PB{\|F}'s minimum generating family.

In practice, the algorithm tends to be interesting only when $m$ and $n$ are
roughly equal. If $n$ is large compared to~$m$, we can remove any vertices of
degree zero; such vertices aren't the endpoint of any interval. If $m$ is
large compared to~$n$, we can almost always find $n$ irredundant intervals by
inspection. The running time in general is readily seen to be $O(mn+n^2)$.

\Y\B\4\D$\\{panic}(\|k)$ \6
${}\{{}$\1\6
${}\\{fprintf}(\\{stderr},\39\.{"Oops,\ we're\ out\ of\ }\)\.{memory!\ (Case\ %
\%d)\\n"},\39\|k);{}$\6
\&{return} \|k;\6
\4${}\}{}$\2\par
\Y\B\8\#\&{include} \.{"gb\_graph.h"}\C{ the GraphBase data structures }\6
\8\#\&{include} \.{"gb\_rand.h"}\C{ the \PB{\\{random\_graph}} generator }\6
\8\#\&{include} \.{"gb\_save.h"}\C{ the \PB{\\{restore\_graph}} generator }\7
\ATH\7
\&{Graph} ${}{*}\|F{}$;\C{ the graph that defines intervals }\6
\&{Graph} ${}{*}\|G{}$;\C{ a graph of intervals that generate \PB{\|F} }\7
\X33:Subroutines\X\7
${}\\{main}(\\{argc},\39\\{argv}){}$\1\1\6
\&{int} \\{argc};\C{ the number of command-line arguments, plus 1 }\6
\&{char} ${}{*}\\{argv}[\,]{}$;\C{ the command-line arguments }\2\2\6
${}\{{}$\1\6
\&{register} \&{Vertex} ${}{*}\|t,{}$ ${}{*}\|u,{}$ ${}{*}\|v,{}$ ${}{*}\|w,{}$
${}{*}\|x{}$;\C{ current vertices of interest }\6
\&{register} \&{Arc} ${}{*}\|a,{}$ ${}{*}\|b,{}$ ${}{*}\|c{}$;\C{ current arcs
of interest }\7
\X24:Scan the command-line options and generate \PB{\|F}\X;\6
\X26:Compute \PB{\|G} and \PB{\|S} by the Franzblau/Kleitman algorithm\X;\6
\&{if} (\\{gb\_trouble\_code})\1\5
\\{panic}(\T{1});\2\6
\X32:Construct an irredundant subfamily of \PB{\|F} with the cardinality of %
\PB{\|G}\X;\6
\X36:Print the results\X;\6
\&{return} \T{0};\C{ this is the normal exit }\6
\4${}\}{}$\2\par
\fi

\M{24}\B\X24:Scan the command-line options and generate \PB{\|F}\X${}\E{}$\6
${}\{{}$\1\6
\&{int} \|m${}\K\T{0},{}$ \|n${}\K\T{0},{}$ \\{seed}${}\K\T{0};{}$\7
\&{if} ${}(\\{argc}\G\T{3}\W\\{sscanf}(\\{argv}[\T{1}],\39\.{"\%d"},\39{\AND}%
\|m)\E\T{1}\W\\{sscanf}(\\{argv}[\T{2}],\39\.{"\%d"},\39{\AND}\|n)\E\T{1}){}$\5
${}\{{}$\1\6
\&{if} ${}(\\{argc}>\T{3}){}$\1\5
${}\\{sscanf}(\\{argv}[\T{3}],\39\.{"\%d"},\39{\AND}\\{seed});{}$\2\6
\&{if} ${}(\|m<\T{0}){}$\5
${}\{{}$\1\6
${}\|m\K{-}\|m{}$;\C{ we assume the user meant to negate \PB{\|n} instead of %
\PB{\|m} }\6
${}\|n\K{-}\|n;{}$\6
\4${}\}{}$\2\6
\&{if} ${}(\|n\G\T{0}){}$\1\5
${}\|F\K\\{random\_graph}(\|n+\T{1},\39\|m,\39\T{0},\39\T{0},\39\T{0},\39\T{0},%
\39\T{0},\39\T{0},\39\T{0},\39\\{seed});{}$\2\6
\&{else}\5
${}\{{}$\1\6
${}\|G\K\\{random\_graph}({-}\|n+\T{1},\39\|m,\39\T{0},\39\T{0},\39\T{0},\39%
\T{0},\39\T{0},\39\T{0},\39\T{0},\39\\{seed});{}$\6
\X25:Set \PB{\|F} to the mirror image of \PB{\|G}\X;\6
\\{gb\_recycle}(\|G);\6
\4${}\}{}$\2\6
\4${}\}{}$\2\6
\&{else} \&{if} ${}(\\{argc}\E\T{2}\W\\{strncmp}(\\{argv}[\T{1}],\39\.{"-g"},%
\39\T{2})\E\T{0}){}$\1\5
${}\|F\K\\{restore\_graph}(\\{argv}[\T{1}]+\T{2});{}$\2\6
\&{else}\5
${}\{{}$\1\6
${}\\{fprintf}(\\{stderr},\39\.{"Usage:\ \%s\ m\ n\ [seed}\)\.{]\ \ |\ \ \%s\
-gfoo.gb\\n"},\39\\{argv}[\T{0}],\39\\{argv}[\T{0}]);{}$\6
\&{return} \T{2};\6
\4${}\}{}$\2\6
\&{if} ${}(\R\|F){}$\5
${}\{{}$\1\6
${}\\{fprintf}(\\{stderr},\39\.{"Sorry,\ can't\ create}\)\.{\ the\ graph!\
(error\ c}\)\.{ode\ \%ld)\\n"},\39\\{panic\_code});{}$\6
\&{return} \T{3};\6
\4${}\}{}$\2\6
\4${}\}{}$\2\6
${}\\{printf}(\.{"Applying\ Franzblau/}\)\.{Kleitman\ to\ \%s:\\n"},\39\|F\MG%
\\{id}){}$;\par
\U23.\fi

\M{25}\B\X25:Set \PB{\|F} to the mirror image of \PB{\|G}\X${}\E{}$\6
$\|F\K\\{gb\_new\_graph}(\|G\MG\|n);{}$\6
\&{if} ${}(\R\|F){}$\1\5
\\{panic}(\T{4});\2\6
${}\\{make\_compound\_id}(\|F,\39\.{"reflect("},\39\|G,\39\.{")"});{}$\6
\&{for} ${}(\|v\K\|G\MG\\{vertices},\39\|u\K\|F\MG\\{vertices}+\|F\MG\|n-%
\T{1};{}$ ${}\|v<\|G\MG\\{vertices}+\|G\MG\|n;{}$ ${}\|u\MM,\39\|v\PP){}$\5
${}\{{}$\1\6
${}\|v\MG\\{clone}\K\|u;{}$\6
${}\|u\MG\\{name}\K\\{gb\_save\_string}(\|v\MG\\{name});{}$\6
\4${}\}{}$\2\6
\&{for} ${}(\|v\K\|G\MG\\{vertices};{}$ ${}\|v<\|G\MG\\{vertices}+\|G\MG\|n;{}$
${}\|v\PP){}$\1\6
\&{for} ${}(\|a\K\|v\MG\\{arcs};{}$ \|a; ${}\|a\K\|a\MG\\{next}){}$\1\5
${}\\{gb\_new\_arc}(\|v\MG\\{clone},\39\|a\MG\\{tip}\MG\\{clone},\39\T{1}){}$;%
\2\2\par
\U24.\fi

\N{1}{26}The main algorithm. We follow the outline of \S\outline.

\Y\B\4\X26:Compute \PB{\|G} and \PB{\|S} by the Franzblau/Kleitman algorithm%
\X${}\E{}$\6
\X27:Make \PB{\|G} a copy of \PB{\|F}, retaining only leftward arcs\X;\6
\&{for} ${}(\|v\K\|G\MG\\{vertices}+\T{2};{}$ ${}\|v<\|G\MG\\{vertices}+\|G\MG%
\|n;{}$ ${}\|v\PP){}$\1\5
\X28:Reduce all minimal bad intervals with right endpoint \PB{\|v} and record
them in \PB{\|S}\X;\2\par
\U23.\fi

\M{27}The algorithm doesn't need pointers from the left endpoint of an interval
to the right endpoint; leftward pointers are sufficient. (This observation
makes the reduction procedure faster.)

While copying \PB{\|F}, we remove its rightward arcs, and we assign
length~1 to all its leftward arcs. Later on, we will represent intervals
of \PB{\|S} by recording them in \PB{\|F} as leftward arcs of length~\PB{${-}%
\T{1}$}.

We also clear two utility fields of \PB{\|F}'s vertices, since they will
be used by the algorithm later. (They might be nonzero, if \PB{\|F} was
supplied with the \.{-g} option.)

\Y\B\4\D$\\{clone}$ \5
$\|u.{}$\|V\C{ the vertex in \PB{\|F} that matches a vertex in \PB{\|G},
         or vice versa }\par
\Y\B\4\X27:Make \PB{\|G} a copy of \PB{\|F}, retaining only leftward arcs\X${}%
\E{}$\6
$\\{switch\_to\_graph}(\NULL){}$;\C{ prepare to return to graph \PB{\|F} later
}\6
${}\|G\K\\{gb\_new\_graph}(\|F\MG\|n){}$;\C{ a graph with nameless vertices and
no arcs }\6
\&{if} ${}(\R\|G){}$\1\5
\\{panic}(\T{5});\2\6
\&{for} ${}(\|u\K\|F\MG\\{vertices},\39\|v\K\|G\MG\\{vertices};{}$ ${}\|u<\|F%
\MG\\{vertices}+\|F\MG\|n;{}$ ${}\|u\PP,\39\|v\PP){}$\5
${}\{{}$\1\6
${}\|u\MG\\{clone}\K\|v{}$;\5
${}\|v\MG\\{clone}\K\|u;{}$\6
${}\|u\MG\|y.\|I\K\|u\MG\|z.\|I\K\T{0};{}$\6
${}\|v\MG\\{name}\K\\{gb\_save\_string}(\|u\MG\\{name});{}$\6
\&{for} ${}(\|a\K\|u\MG\\{arcs},\39\|b\K\NULL;{}$ \|a; ${}\|a\K\|a\MG%
\\{next}){}$\1\6
\&{if} ${}(\|a\MG\\{tip}\G\|u){}$\5
${}\{{}$\C{ we will remove the non-leftward arc \PB{\|a} }\1\6
\&{if} (\|b)\1\5
${}\|b\MG\\{next}\K\|a\MG\\{next};{}$\2\6
\&{else}\1\5
${}\|u\MG\\{arcs}\K\|a\MG\\{next};{}$\2\6
\4${}\}{}$\5
\2\&{else}\5
${}\{{}$\C{ we will copy the leftward arc \PB{\|a} }\1\6
${}\\{gb\_new\_arc}(\|v,\39\|a\MG\\{tip}\MG\\{clone},\39\T{0}){}$;\C{ the
length in $G$ is 0 }\6
${}\|a\MG\\{len}\K\T{1}{}$;\C{ but in $F$ the length is 1 }\6
${}\|b\K\|a{}$;\C{ \PB{\|b} points to the last non-removed arc }\6
\4${}\}{}$\2\2\6
\4${}\}{}$\2\6
\&{if} (\\{gb\_trouble\_code})\1\5
\\{panic}(\T{6});\2\6
\\{switch\_to\_graph}(\|F);\C{ now we can add arcs to \PB{\|F} again }\par
\U26.\fi

\M{28}Here's the most interesting part of the program, algorithmwise.
Given a vertex~\PB{\|v}, we
want to find the largest \PB{\|u}, if any, such that $\[u,v)$ is bad for the
intervals in~\PB{\|G}. So we sweep through the intervals $\[u,v)$ from right to
left, decreasing~\PB{\|u} until it reaches a limiting value~\PB{\|t}. Here \PB{%
\|t} is the
least upper bound on a left endpoint that could guarantee double coverage
of all points in $\[u,v)$.

The utility field \PB{$\|x\MG\\{count}$} records the number of intervals with
left
endpoint \PB{\|x} and right endpoint \PB{\|w} in the range \PB{$\|u<\|w\Z\|v$}.
Another utility field
\PB{$\|x\MG\\{link}$} is used to link vertices with nonzero counts together so
that
we can clear the counts to zero again afterwards.

It's easy to see that each iteration of the \PB{\&{while}} loop in this
section takes at most $O(m+n)$ steps. The actual computation time is,
however, usually much faster.

This program is designed to work correctly when $G$ contains more than one arc
from $u$ to~$v$. Duplicate arcs are discarded as a special case of the
reduction procedure.

\Y\B\4\D$\\{count}$ \5
$\|z.{}$\|I\C{ coverage decreases by this much when we pass to the left }\par
\B\4\D$\\{link}$ \5
$\|y.{}$\|V\C{ pointer to a vertex                   whose \PB{\\{count}} field
needs to be zeroed later }\par
\Y\B\4\X28:Reduce all minimal bad intervals with right endpoint \PB{\|v} and
record them in \PB{\|S}\X${}\E{}$\6
${}\{{}$\1\6
\&{while} (\T{1})\5
${}\{{}$\1\6
\&{int} \\{coverage}${}\K\T{0}{}$;\C{ the number of intervals $\subseteq\[t,v)$
that contain \PB{\|u} }\6
\&{int} \\{potential}${}\K\T{0}{}$;\C{ sum of \PB{$\|x\MG\\{count}$} for \PB{$%
\|x<\|t$} }\6
\&{Vertex} ${}{*}\\{cleanup}\K\NULL{}$;\C{ head of the list of vertices with
nonzero count }\7
\&{for} ${}(\|u\K\|v,\39\|t\K\|v-\T{1};{}$ ${}\|u>\|t;{}$ ${}\|u\MM){}$\5
${}\{{}$\1\6
${}\\{coverage}\MRL{-{\K}}\|u\MG\\{count}{}$;\C{ we prepare to decrease \PB{%
\|u} }\6
\X29:Update the counts for all intervals ending at \PB{\|u}\X;\6
\&{if} ${}(\\{coverage}+\\{potential}<\T{2}){}$\5
${}\{{}$\C{ there's no bad interval ending at \PB{\|v} }\1\6
\X30:Clean up all \PB{\\{count}} fields\X;\6
\&{goto} \\{done\_with\_v};\6
\4${}\}{}$\2\6
\&{while} ${}(\\{coverage}<\T{2}){}$\5
${}\{{}$\1\6
${}\|t\MM;{}$\6
${}\\{coverage}\MRL{+{\K}}\|t\MG\\{count};{}$\6
${}\\{potential}\MRL{-{\K}}\|t\MG\\{count};{}$\6
\4${}\}{}$\2\6
\4${}\}{}$\2\6
${}\\{gb\_new\_arc}(\|v\MG\\{clone},\39\|u\MG\\{clone},\39{-}\T{1}){}$;\C{ $%
\[u,v)$ is minimally bad; we record it in \PB{$\|S\K{-}\|F$} }\6
\X31:Replace $G$ by $G\]\int u,v)$\X;\6
\X30:Clean up all \PB{\\{count}} fields\X;\C{ now we'll try again }\6
\4${}\}{}$\2\6
\4\\{done\_with\_v}:\5
;\6
\4${}\}{}$\2\par
\U26.\fi

\M{29}\B\X29:Update the counts for all intervals ending at \PB{\|u}\X${}\E{}$\6
\&{for} ${}(\|a\K\|u\MG\\{arcs};{}$ \|a; ${}\|a\K\|a\MG\\{next}){}$\5
${}\{{}$\1\6
${}\|w\K\|a\MG\\{tip};{}$\6
\&{if} ${}(\|w\MG\\{count}\E\T{0}){}$\5
${}\{{}$\1\6
${}\|w\MG\\{link}\K\\{cleanup};{}$\6
${}\\{cleanup}\K\|w;{}$\6
\4${}\}{}$\2\6
${}\|w\MG\\{count}\PP;{}$\6
\&{if} ${}(\|w\G\|t){}$\1\5
${}\\{coverage}\PP;{}$\2\6
\&{else}\1\5
${}\\{potential}\PP;{}$\2\6
\4${}\}{}$\2\par
\U28.\fi

\M{30}\B\X30:Clean up all \PB{\\{count}} fields\X${}\E{}$\6
\&{for} ${}(\|w\K\\{cleanup};{}$ \|w; ${}\|w\K\|w\MG\\{link}){}$\1\5
${}\|w\MG\\{count}\K\T{0}{}$;\2\par
\U28.\fi

\M{31}The reduction process is kind of cute too.

\Y\B\4\X31:Replace $G$ by $G\]\int u,v)$\X${}\E{}$\6
\&{for} ${}(\|a\K\|v\MG\\{arcs},\39\|c\K\NULL,\39\|w\K\|v;{}$ \|a; ${}\|c\K\|a,%
\39\|a\K\|a\MG\\{next}){}$\1\6
\&{if} ${}(\|a\MG\\{tip}\G\|u\W\|a\MG\\{tip}<\|w){}$\1\5
${}\|w\K\|a\MG\\{tip},\39\|b\K\|c{}$;\C{ now $\[w,v)$ is the longest interval
from \PB{\|v} inside $\[u,v)$;    we'll remove it }\2\2\6
\&{if} (\|b)\1\5
${}\|b\MG\\{next}\K\|b\MG\\{next}\MG\\{next};{}$\2\6
\&{else}\1\5
${}\|v\MG\\{arcs}\K\|v\MG\\{arcs}\MG\\{next}{}$;\6
\,\C{ the remaining job is to shorten the other maximal arcs in $\[u,v)$ }\2\6
\&{for} ${}(\|t\K\|v-\T{1};{}$ ${}\|w>\|u;{}$ ${}\|t\MM){}$\5
${}\{{}$\1\6
\&{for} ${}(\|a\K\|t\MG\\{arcs},\39\|x\K\|w;{}$ \|a; ${}\|a\K\|a\MG\\{next}){}$%
\1\6
\&{if} ${}(\|a\MG\\{tip}\G\|u\W\|a\MG\\{tip}<\|x){}$\1\5
${}\|x\K\|a\MG\\{tip},\39\|b\K\|a;{}$\2\2\6
\&{if} ${}(\|x<\|w){}$\1\5
${}\|b\MG\\{tip}\K\|w,\39\|w\K\|x{}$;\C{ $\[x,t)$ is the longest interval from %
\PB{\|t} }\2\6
\4${}\}{}$\2\par
\U28.\fi

\N{1}{32}The d\'enouement. Now we build a binary tree in the original graph %
\PB{\|F},
by filling in some of the utility fields of \PB{\|F}'s vertices.
If a node in the tree is labeled with \PB{\|x} and with the interval $[u,v)$,
we represent it by \PB{$\|x\MG\\{left}\K\|u$} and \PB{$\|x\MG\\{right}\K\|v$};
the subtrees of this
node are \PB{$\|x\MG\\{llink}$} and \PB{$\|x\MG\\{rlink}$}. The root of the
whole tree is
\PB{$\|F\MG\\{root}$}.

The \PB{\\{rlink}} field happens to be the same as the \PB{\\{count}} field,
but this
is no problem because the \PB{\\{rlink}} is never changed or examined until
after
the \PB{\\{count}} has been reset to zero for the last time.

\Y\B\4\D$\\{left}$ \5
$\|x.{}$\|V\C{ left endpoint of interval labeling this node }\par
\B\4\D$\\{right}$ \5
$\|w.{}$\|V\C{ right endpoint of interval labeling this node }\par
\B\4\D$\\{llink}$ \5
$\|v.{}$\|V\C{ left subtree of this node }\par
\B\4\D$\\{rlink}$ \5
$\|z.{}$\|V\C{ right subtree of this node }\par
\B\4\D$\\{root}$ \5
$\\{uu}.{}$\|V\C{ root node of the binary tree for this graph }\par
\Y\B\4\X32:Construct an irredundant subfamily of \PB{\|F} with the cardinality
of \PB{\|G}\X${}\E{}$\6
$\|F\MG\\{root}\K\\{make\_tree}(\|F\MG\\{vertices},\39\|F\MG\\{vertices}+\|F\MG%
\|n-\T{1}){}$;\par
\U23.\fi

\M{33}With a little care we could maintain a stack within \PB{\|F} itself, but
it's easier to use recursion in \CEE/. Let's just hope the system
programmers have given us a large enough runtime stack to work with.

This subroutine is based on the trick explained in \S\trick.

\Y\B\4\X33:Subroutines\X${}\E{}$\6
\&{Vertex} ${}{*}\\{make\_tree}(\|t,\39\|w){}$\1\1\6
\&{Vertex} ${}{*}\|t,{}$ ${}{*}\|w;\2\2{}$\6
${}\{{}$\1\6
\&{register} \&{Vertex} ${}{*}\|u,{}$ ${}{*}\|v,{}$ ${}{*}\|x;{}$\6
\&{register} \&{Arc} ${}{*}\|a;{}$\7
\X34:Find a vertex \PB{\|x} with $N_x\,F\:\int t,w)-N_x\,S\:\int t,w)=1$\X;\6
\&{if} ${}(\R\|x){}$\1\5
\&{return} ${}\NULL{}$;\C{ $F\:\[t,w)$ is empty }\2\6
\X35:Find an interval $\int u,v)$ such that $x\in\int u,v)\seq\int t,w)$\X;\6
${}\|x\MG\\{left}\K\|u;{}$\6
${}\|x\MG\\{right}\K\|v;{}$\6
${}\|x\MG\\{llink}\K\\{make\_tree}(\|t,\39\|x);{}$\6
${}\|x\MG\\{rlink}\K\\{make\_tree}(\|x+\T{1},\39\|w);{}$\6
\&{return} \|x;\6
\4${}\}{}$\2\par
\A37.
\U23.\fi

\M{34}A subtle bug is avoided here when we realize that a vertex might
already be in the cleanup list when its \PB{\\{count}} is zero.

\Y\B\4\X34:Find a vertex \PB{\|x} with $N_x\,F\:\int t,w)-N_x\,S\:\int t,w)=1$%
\X${}\E{}$\6
${}\{{}$\1\6
\&{register} \&{int} \\{coverage}${}\K\T{0}{}$;\C{ coverage in \PB{\|F} minus %
\PB{\|S} }\6
\&{Vertex} ${}{*}\\{cleanup}\K\|w+\T{1}{}$;\C{ \PB{$\|w+\T{1}$} is a sentinel
value }\7
\&{for} ${}(\|v\K\|w,\39\|x\K\NULL;{}$ ${}\|v>\|t;{}$ ${}\|v\MM){}$\5
${}\{{}$\1\6
${}\\{coverage}\MRL{-{\K}}\|v\MG\\{count}{}$;\C{ now \PB{\\{coverage}} refers
to $N_{v-1}$ }\6
\&{for} ${}(\|a\K\|v\MG\\{arcs};{}$ \|a; ${}\|a\K\|a\MG\\{next}){}$\5
${}\{{}$\1\6
${}\|u\K\|a\MG\\{tip};{}$\6
\&{if} ${}(\|u\G\|t){}$\5
${}\{{}$\1\6
\&{if} ${}(\|u\MG\\{link}\E\NULL){}$\1\5
${}\|u\MG\\{link}\K\\{cleanup},\39\\{cleanup}\K\|u;{}$\2\6
${}\|u\MG\\{count}\MRL{+{\K}}\|a\MG\\{len}{}$;\C{ the length is $+1$ for $F$,
$-1$ for $S$ }\6
${}\\{coverage}\MRL{+{\K}}\|a\MG\\{len};{}$\6
\4${}\}{}$\2\6
\4${}\}{}$\2\6
\&{if} ${}(\\{coverage}\E\T{1}){}$\5
${}\{{}$\1\6
${}\|x\K\|v-\T{1}{}$;\5
\&{break};\6
\4${}\}{}$\2\6
\4${}\}{}$\2\6
\&{if} ${}(\R\|x\W\\{cleanup}\Z\|w){}$\1\5
${}\\{fprintf}(\\{stderr},\39\.{"This\ can't\ happen!\\}\)\.{n"});{}$\2\6
\&{while} ${}(\\{cleanup}\Z\|w){}$\5
${}\{{}$\1\6
${}\|v\K\\{cleanup}\MG\\{link};{}$\6
${}\\{cleanup}\MG\\{count}\K\T{0};{}$\6
${}\\{cleanup}\MG\\{link}\K\NULL;{}$\6
${}\\{cleanup}\K\|v;{}$\6
\4${}\}{}$\2\6
\4${}\}{}$\2\par
\U33.\fi

\M{35}\B\X35:Find an interval $\int u,v)$ such that $x\in\int u,v)\seq\int
t,w)$\X${}\E{}$\6
\&{for} ${}(\|v\K\|w;{}$ ${}\|v>\|x;{}$ ${}\|v\MM){}$\5
${}\{{}$\1\6
\&{for} ${}(\|a\K\|v\MG\\{arcs};{}$ \|a; ${}\|a\K\|a\MG\\{next}){}$\1\6
\&{if} ${}(\|a\MG\\{len}>\T{0}){}$\5
${}\{{}$\1\6
${}\|u\K\|a\MG\\{tip};{}$\6
\&{if} ${}(\|u\Z\|x\W\|u\G\|t){}$\1\5
\&{goto} \\{done};\2\6
\4${}\}{}$\2\2\6
\4${}\}{}$\2\6
\\{done}:\par
\U33.\fi

\M{36}\B\X36:Print the results\X${}\E{}$\6
\\{printf}(\.{"Minimum\ generating\ }\)\.{family:"});\6
\&{for} ${}(\|v\K\|G\MG\\{vertices}+\T{1};{}$ ${}\|v<\|G\MG\\{vertices}+\|G\MG%
\|n;{}$ ${}\|v\PP){}$\1\6
\&{for} ${}(\|a\K\|v\MG\\{arcs};{}$ \|a; ${}\|a\K\|a\MG\\{next}){}$\1\5
${}\\{printf}(\.{"\ [\%s..\%s)"},\39\|a\MG\\{tip}\MG\\{name},\39\|v\MG%
\\{name});{}$\2\2\6
\\{printf}(\.{"\\nMaximum\ irredunda}\)\.{nt\ family:"});\6
${}\\{postorder\_print}(\|F\MG\\{root});{}$\6
\\{printf}(\.{"\\n"});\par
\U23.\fi

\M{37}There's just one subroutine to go. This is textbook stuff.

\Y\B\4\X33:Subroutines\X${}\mathrel+\E{}$\6
\&{void} \\{postorder\_print}(\|x)\1\1\6
\&{Vertex} ${}{*}\|x;\2\2{}$\6
${}\{{}$\1\6
\&{if} (\|x)\5
${}\{{}$\1\6
${}\\{postorder\_print}(\|x\MG\\{llink});{}$\6
${}\\{postorder\_print}(\|x\MG\\{rlink});{}$\6
${}\\{printf}(\.{"\ \%s[\%s..\%s)"},\39\|x\MG\\{name},\39\|x\MG\\{left}\MG%
\\{name},\39\|x\MG\\{right}\MG\\{name});{}$\6
\4${}\}{}$\2\6
\4${}\}{}$\2\par
\fi

\N{1}{38}Comments and extensions.
The program just presented incorporates several refinements to the
implementation sketched by Franzblau and Kleitman in~[\FK], and the author
hopes that readers will enjoy finding them in the code. The Stanford GraphBase
provides convenient data structures, by means of which it was possible to
make the program short and sweet. However, most of the
key ideas (except for the \PB{\\{make\_tree}} procedure) can be found in~[\FK].

\fi

\M{39}Lubiw [\Lu] discovered that the algorithm of Franzblau and Kleitman can
be
generalized so that it finds optimum irredundant subfamilies and
generating families in an appropriate sense when the points of the underlying
line have been given arbitrary nonnegative {\it weights}. It should be
interesting and instructive to extend the program above so that it
handles this more general problem.

\fi

\M{40}The introductory remarks in \S\intro\ mention the recent breakthrough by
Frank and Jord\'an~[\FJ], who showed (among many other things) that
Gy\H{o}ri's theorem can be extended to intervals on a circle as well as a
line. Such a generalization was surprising because the size of a minimum
generating set might be strictly larger than the size of a maximum irredundant
subfamily of cyclic intervals. For example, the $n$ intervals $\[k,k+2)$ for
$0\le k<n$ on
the ring of integers mod~$n$ are obviously redundant; if we leave out any one
of them, the remaining $n-1$ intervals
will cover all $n$ points. However, these cyclic
intervals cannot be generated by fewer than $n$ subintervals: No $n-1$
subintervals of length~1 will do the job, and if $\[k,k+2)$ is one of the
generators the remaining $n-1$ intervals require $n-1$ further generators
because they are irredundant.

\callthis\firstext
Gy\H{o}ri's minimax principle is restored, however, if we change the
definition of irredundant families.
We can say that $F$ is irredundant if each $f\in F$ has a distinguished
element $f_*\in f$, such that whenever $f$ and~$f'$ are distinct sets of~$F$
we have either $f_*\notin f'$ or $f'_*\notin f$. If $F$ is irredundant in this
sense, and if $G$~generates~$F$, it is not difficult to prove that
$\vert F\vert\le\vert G\vert$: There is a $g_f\in G$ for each $f\in F$, with
the property that $f_*\in g_f\subseteq f$; our new definition guarantees that
$g_f\ne g_{f'}$ when $f\ne f'$.

According to this new definition, the intervals $\[k,k+d)$ for $0\le k<n$,
modulo~$n$, are irredundant whenever $n>2(d-1)$, because we can let
$\[k,k+d)_*=k$. Frank and Jord\'an showed that if $F$ is any family of
intervals modulo~$n$ with the property that each intersection $f\cap f'$ of
two of its
members is either empty or a single interval, then the size of $F$'s smallest
generating family is the size of its largest irredundant subfamily under the
new definition.

\fi

\M{41}For intervals on a line, Gy\H{o}ri [\Gy] had already observed that both
definitions of irredundancy are equivalent. Suppose a system of
representatives $f_*\in f$ is given for all $f$ in some family~$F$ of
intervals on a line, such that $f\ne f'$ implies $f_*\notin f'$ or $f'_*\notin
f$. If we cannot arrange those intervals in a sequence $f\(1)$, $f\(2)$,
\dots, $f\(n)$ such that $f_*\(j)\notin f\(1)\cup\cdots\cup f\(j-1)$
for $1<j\le n$, there must be some cycle of intervals such that
$f_*\(1)\in f\(2)$, $f_*\(2)\in f\(3)$, \dots, $f_*\(m)\in f\(1)$,
where $m>2$. Consider the shortest such cycle, and suppose
$f_*\(1)=\min_{j=1}^m f_*\(j)$. We cannot have $f\(k-1)<f_*\(k)$ for
$1<k\le m$, because $f_*\(m)\in f\(1)$; let $k>1$ be minimum such that
$f\(k-1)$ is not strictly less than~$f_*\(k)$. Then
$f\(k-1)$ must be strictly greater than~$f_*\(k)$, and we have
$f_*\(1)<f_*\(k)<f_*\(k-1)$. There is some $j$ with $1<j<k$ and
$f_*\(j-1)<f_*\(k)<f_*\(j)$; since $f_*\(j-1)\in f\(j)$ we have
$f_*\(k)\in f\(j)$, a shorter cycle. This contradiction shows that no
cycles exist.

\fi

\M{42}\callthis\secondext
Frank and Jord\'an gave another criterion for irredundancy that works
also for general families of intervals on a circle when large intervals might
wrap around so that their intersection$f\cap f'$ consists of two disjoint
intervals. In such cases they allow $\{f_*,f'_*\}\subseteq f\cap f'$, but only
if $f_*$ and $f'_*$ lie in different components of $f\cap f'$. For example,
the intervals $\[k,k+d)$ for $0\le k<n$, modulo~$n$, are irredundant by this
definition for all $n>d$. Once again the minimax theorem for
generating families and irredundant subfamilies remains valid, in this
extended sense.

\fi

\M{43}The algorithms presented by Frank and Jord\'an [\FJ] for such problems
require linear programming as a subroutine. Therefore it would be extremely
interesting to find a purely combinatorial procedure, analogous to the
algorithm of Franzblau and Kleitman, either for the wrap-restricted situation
of \S\firstext\ or for the more general setup of \S\secondext.

\fi

\N{1}{44}References.

\newdimen\refindent \refindent=27pt
\def\ref[#1] {\smallskip\noindent\hangindent\refindent
\hbox to\refindent{\hfill[#1]\enspace}}

\ref[\CR] Joseph C. Culberson and Robert A. Reckhow, ``Covering polygons
is hard,'' {\sl Journal of Algorithms\/ \bf17} (1994), 2--44.

\ref[\FJ] Andr\'as Frank and Tibor Jord\'an, ``Minimal edge-coverings of pairs
of sets,'' {\sl Journal of Combinatorial Theory\/ \bf B65} (1995), 73--110.

\ref[\FK] Deborah S. Franzblau and Daniel J. Kleitman, ``An algorithm for
constructing polygons with rectangles,'' {\sl Information and Control\/
\bf63} (1984), 164--189.

\ref[\Gy] Ervin Gy\H{o}ri, ``A minimax theorem on intervals,'' {\sl Journal
of Combinatorial Theory\/ \bf B37} (1984), 1--9.

\ref[\KB] Donald E. Knuth and P. B. Bendix, ``Simple word problems in
universal algebras,'' in {\sl Computational Problems in Abstract
Algebra}, edited by J. Leech (Oxford:  Pergamon, 1970), 263--297.
Reprinted in {\sl Automation of Reasoning}, edited by J\"org~H.
Siekmann and Graham Wrightson, {\bf2} (Springer, 1983), 342--376.

\ref[\vone] Donald E. Knuth, {\sl Fundamental Algorithms}, Volume~1 of
{\sl The Art of Computer Programming\/} (Reading, Massachusetts:
Addison-Wesley, 1968), $\rm xxii + 634$~pp.

\ref[\vthree] Donald E. Knuth, {\sl Sorting and Searching}, Volume~3 of
{\sl The Art of Computer Programming\/} (Reading, Massachusetts:
Addison-Wesley, 1973), $\rm xii + 722$~pp.\ + foldout illustration.

\ref[\sgb] Donald E. Knuth, {\sl The Stanford GraphBase\/}:\ A Platform
for Combinatorial Computing (New York: ACM Press, 1993), $\rm viii+576$~pp.
Available via anonymous ftp from \.{labrea.stanford.edu} in
directory \.{\char`\~ftp/pub/sgb}.

\ref[\cweb] Donald E. Knuth and Silvio Levy, {\sl The {\slttten CWEB} System
of Structured Documentation\/} (Reading, Massachusetts: Addison-Wesley,
1994).

\ref[\Lu] Anna Lubiw, ``A weighted min-max relation for intervals,''
{\sl Journal of Combinatorial Theory\/ \bf B53} (1991), 173--194.

\fi

\M{45}The author thanks the referees of this note for many valuable suggestions
that greatly improved the presentation.

\fi

\N{1}{46}Index.

\fi

\inx
\fin
\con